\documentclass{article}
\usepackage{indentfirst, amsfonts, amsmath, amsthm,amssymb,amscd}
\usepackage[dvips]{graphicx}

\begin{document}

\title{\bf Compact 4-manifolds with harmonic Weyl tensor and nonnegative biorthogonal curvature}

\author{ Ezio Araujo Costa}
\date{}
\maketitle
\begin{center}{\bf Abstract}
\end{center}
 {\bf  In this note we classify  compact 4-manifolds with harmonic Weyl tensor and nonnegative biorthogonal curvature }
 \\
 \\
{\bf M. S. C. (2000)}: 53C25, 53C24.
\\
\\
Let $(M^4, g)$ be a Riemannian  4-manifold,
$s$  the scalar curvature of this metric $g$ and $K$ its sectional curvature. For each $p\in M^4$, let  $P$ a 2-plain in the tangent space $T_pM^4$ and $K^\perp(P)$ its orthogonal complement in $T_pM^4$. The {\it biorthogonal (sectional) curvature } relative to $P$ (in $p\in M^4$) is the average
$$K^\perp (P) = \frac{K(P) + K(P^\perp)}{2}\eqno [1.1]$$
 Then we have the following function in $M^4$:

 $$K_1^\perp(p)= inf \{K^\perp (P);  P \subset T_pM^4\},\eqno [1.2]$$

If $(M^4, g)$ be an oriented Riemannian  4-manifold, the Weyl tensor has the decomposition $W = W^+\bigoplus W^-$, where $W^\pm : \Lambda^\pm
\longrightarrow \Lambda^\pm$ are self-adjoint with free traces and are called of the self-dual and anti-self-dual parts of $W$, respectively. Let $w_1^\pm \leq w_2^\pm \leq w_3^\pm$ be the eigenvalues of $W^\pm$, respectively.
\\
As was proved in [3] we have the following relation:
$$K_1^\perp - \frac{s}{12} = \frac{w_1^+ + w_1^-}{2}, \eqno [1.3]$$

A Riemannian manifold $(M, g)$ is said to have {\it harmonic curvature} if its Levi-Civita connection $\nabla$
in the tangent bundle $TM$ satisfies $d^*R = 0$, where $R$ is the curvature tensor of $(M, g)$ and $d$ is the operator of exterior differentiation. If $M$ is compact, this just means that $\nabla$ is a critical point for the Yang-Milles functional $YM(\nabla) = \frac{1}{2}\int_{M}\mid R\mid^2$, where $R$ is the curvature of connection $\nabla$. Einstein metrics has harmonic curvature and metrics with harmonic curvature has constant scalar curvature .

Let $(M^4 ,g)$ be an oriented Riemannian 4-manifold with self-dual Weyl tensor $W^+$ and anti-self-dual Weyl tensor $W^-$.
Viewing $W^\pm$ as a tensor of type (0,4), we say that $W^+$ ($W^-$) is harmonic if $\delta W^+ = 0$ ($\delta W^+ = 0$, respec.), where $\delta$ is the formal divergence defined for any tensor $T$ of type (0,4) by $$\delta T(X_1,X_2,X_3) = -trace_{g}\{(Y,Z)\mapsto\nabla_{Y}T(Z,X_1,X_2,X_3)\},$$ where $g$ is the metric of $M$. $(M,g)$ has harmonic Weyl tensor $W$ if $W^\pm$ are harmonics.
Einstein metrics and metrics with harmonic curvature are analytical metrics and has harmonic tensor Weyl.

A Riemannian manifold $(M, g)$ is said to have {\it harmonic curvature} if its Levi-Civita connection $\nabla$
in the tangent bundle $TM$ satisfies $d^*R = 0$, where $R$ is the curvature tensor of $(M, g)$ and $d$ is the operator of exterior differentiation. If $M$ is compact, this just means that $\nabla$ is a critical point for the Yang-Milles functional $YM(\nabla) = \frac{1}{2}\int_{M}\mid R\mid^2$, where $R$ is the curvature of connection $\nabla$. Notice that Einstein metrics has harmonic curvature and metrics with harmonic curvature has constant scalar curvature.
Let $(M^4,g)$ be an oriented Riemannian 4-manifold with self-dual Weyl tensor $W^+$ and anti-self-dual Weyl tensor $W^-$.
Viewing $W^\pm$ as a tensor of type (0,4), we say that $W^+$ ($W^-$) is harmonic if $\delta W^+ = 0$ ($\delta W^+ = 0$, respec.), where $\delta$ is the formal divergence defined for any tensor $T$ of type (0,4) by $$\delta T(X_1,X_2,X_3) = -trace_{g}\{(Y,Z)\mapsto\nabla_{Y}T(Z,X_1,X_2,X_3)\},$$ where $g$ is the metric of $M$. $(M^4,g)$ has harmonic Weyl tensor $W$ if $W^\pm$ are harmonics.
 Einstein metrics and metrics with harmonic curvature are analytical metrics and has harmonic tensor Weyl. In [1], R. Bettiol obtained the topological classification of compact simply connected  4-manifolds with positive biorthogonal  and compacts 4-manifolds with harmonic Weyl tensor and positive biortoghonal  was studied in [3] end [5]. In this note we prove the following
\\
\\
{\bf Theorem 1.1}-{\it Let $(M^4, g$) be a compact oriented Riemannian 4-manifold with harmonic Weyl tensor, analytical metric and nonnegative biorthogonal curvature. Then $(M^4, g$) is an Einstein manifold with nonnegative sectional curvature, $(M^4, g$) is conformally flat or the universal covering of $(M^4, g$) is diffeomorphic to $\mathbb{S}^2 \times \mathbb{S}^2$ and $g$ is a product of metrics with nonnegative sectional curvature.}

\begin{center}{\bf Proof}
\end{center}
By the proof of Theorem 6 in [3], we have $$\mid W^+ \mid + \mid W^-\mid \leq \sqrt{6}[s/6 -2K_1^\perp] \leq s/\sqrt{6},$$ 
where $s$ is the scalar curvature of $(M^4, g)$.
\\
Let $A =\{ p\in M^4; Ric(p) \neq s(p)/4\}$ and $B = \{p \in M^4; \mid W^+ \mid = \mid W^-\mid\}$. 
\newpage
Note that $A \subset B$. If $A$ is empty then $(M^4, g)$ is an Einstein manifold.
Otherwise we deduce that $B = M^4$ and using the Weitzenbock formula for $W^\pm$, is easy see that $W^\pm =0$ or $\mid W^+ \mid^2 = \mid W^-\mid^2 = s^2/24$, $w_1\pm = w_2^\pm$ and so $(M^4,g)$ has nonnegative isotropic curvature. In this case we can uses the theorem B in [2] (see also Prop. 3 in [4])

\end{document}